\documentclass[10]{amsart}
\usepackage{amssymb, amsmath, amsfonts, amsthm, graphics,hyperref, footmisc, blindtext,latexsym,stmaryrd, enumerate}
\usepackage[hmargin=1 in, vmargin = 1.5 in]{geometry}
\usepackage{hyperref}
\usepackage{color} 
\usepackage[usenames,dvipsnames]{xcolor} 
\definecolor{darkblue}{rgb}{0,0,1}
\usepackage{hyperref} 
\hypersetup{colorlinks,breaklinks,linkcolor=red,urlcolor=darkblue,
anchorcolor=darkblue,citecolor=green}

\usepackage[all]{xy}
\usepackage{tikz-cd}
\usepackage{blindtext}

\usepackage{ytableau}
\numberwithin{figure}{section}
\numberwithin{equation}{section}

\newcommand{\into}{\hookrightarrow}

\newcommand{\CC}{\mathbb{C}}

\renewcommand{\SS}{\mathbb{S}}

\newcommand{\ZZ}{\mathbb{Z}}

\newcommand{\cZ}{\mathcal{Z}}

\newcommand{\gl}{\operatorname{GL}}

\newtheorem{theorem}[figure]{Theorem}
\newtheorem{lemma}[figure]{Lemma}
\newtheorem{corollary}[figure]{Corollary}
\newtheorem{proposition}[figure]{Proposition}
\theoremstyle{definition}
\newtheorem{definition}[figure]{Definition}
\newtheorem{notation}[figure]{Notation}
\newtheorem{construction}[figure]{Construction}

\theoremstyle{remark}
\newtheorem{remark}[figure]{Remark}
\newtheorem{note}[figure]{Note}
\newtheorem{example}[figure]{Example}

\title{Composition Tableaux basis for Schur functors and the Pl\"ucker algebra}
\author{Shubhankar Sahai}

\begin{document}

\begin{abstract}
   We show that combinatorial objects called row-strict composition tableaux, introduced by Mason and Remmel in 2014 and closely related to the quasi-symmetric   Schur functions of Haglund–Luoto–Mason–van Willigenburg, form a basis for Schur functors of finite free modules over arbitrary commutative rings.
   
   When the ring is the complex numbers, this produces a new basis for the irreducible polynomial representations of $\operatorname{GL}_n(\CC).$ Moreover, in this case it also produces new basis for the Pl\"ucker algebra, a subalgebra of the polynomial ring over $\CC$ in $n^2$ variables, which is of independent combinatorial and geometric interests.
   
   As an aside we also show that these results hold for other combinatorial objects called reverse row strict tableau.
\end{abstract}.
\maketitle

\section{Introduction}

The aim of this note is to explain a representation theoretic interpretation of the \emph{row-strict composition tableaux} introduced by Mason and Remmel in \cite{mason} (see also the theory of \emph{quasi-symmetric Schur functions} of \cite{Steph}). To state our results, it is useful to first recall the classical theory first.

\subsection{Classical theory} This theory is developed in detail in Chapter $8$ of Fulton's book \cite{fulton}.

Let $R$ be a commutative ring. Given a partition $\lambda=(\lambda_1,\ldots, \lambda_k)$ of a non-negative integer, we denote by $\SS^\lambda(-)$ the \emph{Schur functor labelled by $\lambda$}. These are endo-functors of the category of $R$-modules which specialize to the usual exterior and symmetric powers or $R$-modules for specific $\lambda$.  $\SS^{\lambda}(-)$, of course, depends on the partition $\lambda$ but we keep the choice implicit as there is little danger of confusion.

If $E$ is a finite free $R$-module, then so is $\SS^\lambda(E)$ and so one can regard $\SS^\lambda(-)$ as an endofunctor on the category of finite free $R$-modules. Now for the rest of this section we fix a finite free $R$-module $E.$

With $E$ as above, there are several equivalent ways to give  an explicit construction of $\SS^\lambda(E)$. The one most pertinent to our representation theoretic goal is to build it as a certain sub-modules of the polynomial ring $R[Z]:=R[Z_{1,1},Z_{1,2},\ldots, Z_{n,n}]$ by identifying $\SS^\lambda(E)$ with the $R$-module generated by the products of top justified matrix minors of the generic $n\times n$ matrix, $[Z_{i,j}]$ (i.e. the matrix with entries in the $n\times n$ indeterminates $\{Z_{1,1},Z_{1,2},\ldots, Z_{n,n}\}$ ) labelled by the \textit{semistandard Young's tableaux} of shape $\lambda$ with entries in $[n]:=\{1,\ldots, n\}$.

Now suppose $R=\mathbb{C}$, the field of complex numbers. In this case $E$ is a finite dimensional vector space over $\CC$ and we set $n=\operatorname{dim}_{\CC}(E).$ Now, by functoriality, $\SS^\lambda(E)$ has a canonical action of $\gl(E)$, the general linear group of $E$ and thus $\SS^\lambda(E)$ is canonically a $\gl(E)$ representation. This representation of $\gl(E)$ turns out to be irreducible. Moreover, it is a polynomial representation in the sense that upon choosing coordinates for $E$, the homomorphism $\gl(\CC^n)\to \gl(\SS^\lambda(\CC^n))$ is given by polynomial functions on the algebraic variety defined by $\gl(\CC^n)$.

Moreover any other irreducible polynomial representation of $\gl(E)$ can be constructed as $\SS^{\lambda'}(E)$ as $\lambda'$ ranges over partitions of non-negative integers of length at most $n$.

In this case one say a little bit more about the embedding of the $\SS^\lambda(E)$ inside $\CC[X]$. Note that $\operatorname{End}(E)$, the endomorphisms of $E$, is a $\gl(E)\times \gl(E)$ representation and therefore so is $\operatorname{Sym}^d(\operatorname{End}(E))$, the $d$-th symmetric power of $\operatorname{End}(E)$. Then for each non-negative $d$ there is a functorial decomposition of $\operatorname{Sym}^d(\operatorname{End}(E))$ as $\gl(E)\times \gl(E)$ representations as
\begin{equation}
    \operatorname{Sym}^d(\operatorname{End}(E))\cong \bigoplus_{\lambda \vdash d, l(\lambda)\leq n} \mathbb{S}^\lambda(E^\vee)\otimes_{\CC}\mathbb{S}^\lambda (E)
,\end{equation}
where $\lambda \vdash d$ means that $\lambda$ is a partition of $d$ and $l(\lambda)$ is the length of $\lambda$. This implies the decomposition of the symmetric algebra $\operatorname{Sym}^*(\operatorname{End}(E))$ as 
\begin{equation}
    \operatorname{Sym}^*(\operatorname{End}(E))\cong \bigoplus_{l(\lambda) \leq  n} \mathbb{S}^\lambda(E^\vee)\otimes_{\CC}\mathbb{S}^\lambda (E)
\end{equation}

Choosing coordinates we see that the left hand side can be identified with the ring $\CC[Z]$. One then has a decomposition as $\gl_n(\CC)\times \gl_n(\CC)$ representations,

\begin{equation}
    \CC[Z]\cong \bigoplus_{l(\lambda)\leq n}\mathbb{S}^\lambda((\CC^{n})^\vee)\otimes_{\CC}\mathbb{S}^\lambda(\CC^n). 
\end{equation}

Let $N_{-}$ denote unipotent radical of the Borel subgroup of $\gl_n(\CC)$ i.e. the subgroup consisting of lower triangular matrices with $1$'s on the diagonal. Then the fixed ring $\CC[Z]^{N_{-}\times id_{\gl_n(\CC)}}$ is the sub-ring generated by the $2^n$ minors of the matrix $[Z_{i,j}]$. It is called the \emph{Pl\"ucker algebra} and is denoted $\Pi$.  Using highest weight theory, we have an identification

\begin{equation}\label{Pluckeralgebra}
    \Pi \cong \bigoplus_{l(\lambda)\leq n}\mathbb{S}^\lambda(\CC^n). 
\end{equation}

The above story then translates to the result that the polynomials labelled by the semi-standard Young's tableaux form a $\CC$-basis of the Pl\"ucker algebra.

\subsection{Statement of main results}

We now give a brief overview of our results. We will give more details in Section \ref{prelims}.

Let $\alpha=(\alpha_1,\ldots, \alpha_l)$ be a composition of an integer with parts at most $n$ i.e. $\alpha_i\leq n$ for all $i$. Note that every composition can be rearranged into a unique partition by imposing that the tuple be ordered. We denote this partition by $P(\alpha)$.
Analogous to Young's diagrams, one can associate to each such composition a composition diagram - an array with $\alpha_i$ rows in the $i$th row. Note that the composition diagram of $P(\alpha)$ considered as a composition is the Young's diagram of $P(\alpha)$ considered as a partition.
In \cite{mason}, the authors introduce a filling of such composition diagrams, called the row-strict composition tableaux. Let $\mathcal{RSCT}(\alpha,[n])$ denote the set of all row-strict composition tableaux fillings of $\alpha$ with values in $[n]$. We describe a method of uniquely associating to each element of $\mathcal{RSCT}(\alpha,[n])$ a product of the top justified minors of the matrix $[Z_{i,j}]$. We can now state our main theorem.

\begin{theorem}\label{theorem1}
Let $R$ be a commutative ring and $E$ a finite free $R$-module of rank $n$. Let $\lambda$ be a partition with length at most $n$ and $\mathbb{S}^\lambda(-)$ the associated Schur functor. Let $\lambda^T$ be the transpose of $\lambda$ and consider the set $$\mathcal{RSCT}(\lambda^T,[n]):=\bigsqcup_{P(\alpha)=\lambda^T}\mathcal{RSCT}(\alpha,[n]).$$
Then the products of top justified minors associated to the set $\mathcal{RSCT}(\lambda^T,[n])$\footnote{The reason we have to use $\lambda^T$ instead of $\lambda$ is enforced upon us by the nature of the RSCT fillings.} form a basis for the embedding of $\mathbb{S}^\lambda(E)\into R[Z]$.
\end{theorem}

As an immediate application we obtain a representation theoretic corollary when we specialise to the case of $\CC$,

\begin{corollary}\label{repcor}
Let $R=\CC$ and keeping notation as in Theorem \ref{theorem1} . The top justified minors labelled by the set $\mathcal{RSCT}(\lambda,[n])$ form a basis for the embedding of the $\gl(E)$ representations $\mathbb{S}^\lambda(E)$ inside $\CC[Z].$
\end{corollary}

From (\ref{Pluckeralgebra}) we obtain another corollary of representation theoretic importance,

\begin{corollary}\label{basisplucker}
Let $R=\CC$ and keep notation as in Theorem \ref{theorem1}. Let $\Pi$ be the Pl\"ucker algebra. Then the products of top justified matrix minors labelled by $\mathcal{RSCT}(\lambda,[n])$ form a basis for $\Pi$ as a $\CC$-vector space.
\end{corollary}

In the development of the proof of this result we also show that the \textit{reverse-row strict composition tableaux} of Definition \ref{definitionrsst} satisfies analogous results, which is probably known to experts.

\begin{theorem}\label{theorem2}
Fixing notation as in Theorem \ref{theorem1}, we let $\mathcal{RRST}(\lambda^T,[n])$ denote the set of all reverse-row strict fillings of $\lambda^T$ with values in $[n].$ Then the products of top justified minors associated to the set $\mathcal{RRST}(\lambda^T,[n])$ form a basis for the embedding of $\mathbb{S}^\lambda(E)\into R[Z]$.
\end{theorem}

\subsection{Acknowledgements}

The author is indebted to David Speyer for suggesting this problem in the 2017 summer REU program at the University of Michigan, Ann Arbor. He is also grateful to Visu Makam for discussions concerning the representation theory involved. Lastly he is deeply indebted to the anonymous referee for several corrections and observations which greatly streamlined the presentation of the paper and most importantly simplified the proof of the main theorem.
\section{Preliminaries}\label{prelims}
The aim of this section is to set notation and recall some classical definitions in our context.
In this section we fix a ring $R$ and a finite free $R$-module $E$. We let $n$ be the rank of $E$.
\subsection{Matrix minors of a multi-index}

\begin{construction}
Let $I=(i_1,\ldots, i_k)$ be a multi-index with entries in $[n]$. To such an $I$ one can associate a matrix minor of the matrix of indeterminates $[Z_{i,j}]$ in the following manner. Suppose the cardinality of $I$, denoted $|I|$ is $k\geq 0$. Then we look at the top $k\times k$ justified minor with columns labelled by the entries of $I.$ This matrix minor is called the \textit{Pl\"ucker coordinate} of $I$ and is denoted $\Delta_I$.
\end{construction}

\begin{example}
Let $(1,3,4)$ be a multi-index with values in $[n]$ for $n$ large enough. Then the associated Pl\"ucker coordinate $\Delta_{(1,3,4)}$ is given by the determinant 

$$\textrm{det}\begin{pmatrix}
    Z_{1,1} & Z_{1,3} &Z_{1,4}  \\
    Z_{2,1} & Z_{2,3} &Z_{2,4}\\
    Z_{3,1}&Z_{3,3}& Z_{3,4}
    \end{pmatrix}.$$
\end{example}

\begin{remark}
Note that for a mult-index $I=(i_1,\ldots, i_k)$, the associated Pl\"ucker coordinate might be $0$ incase $i_m=i_n$ when $m\neq n$. Further note that there is an action of $S_k$ on the multi-index $(i_1,\ldots, i_k)$, but this only changes the associated Pl\"ucker coordinate by a sign change. This will not affect the results in this paper. Therefore for concreteness we stipulate that any multi-index $I=(i_1,\ldots, i_k)$ be arranged in weakly increasing order so that $(i_1\leq \ldots\leq i_k)$. 
\end{remark}
\subsection{Semistandard Young Tableaux}

Recall that to any partition $\lambda=(\lambda_1,\ldots, \lambda_k)$ one can associate a Young's diagram, an left justified array of containing $\lambda_i$ cells in the $i$th row.

\begin{example}
Let $\lambda=(5,3,1)$ then the Young's diagram is given by
$$
\ytableausetup{centertableaux}\ytableausetup{nosmalltableaux}
\ydiagram{5,3,1}.$$
\end{example}

In any such Young's diagram $\lambda$ one can fill the boxes with entries from positive integers. This is called a \textit{filling of shape $\lambda$}.

\begin{example}\label{wrongssyt}
A filling of shape $(5,3,1)$ is 
$$
\ytableausetup{nosmalltableaux}
\begin{ytableau}
1 & 7 & 7&2&3 \\
2 & 3  &8\\
5\\
\end{ytableau}$$
\end{example}

\begin{notation}
For record keeping purposes, it will be useful to denote the set of all fillings of fixed shape $\lambda $ with values in the set $[n]$ as $\mathcal{FILL}(\lambda,[n])$.
\end{notation}

\begin{construction} \label{constfilling}
For a filling $J\in \mathcal{FILL}(\lambda,[n])$ we will construct a product of Pl\"ucker monomials in two ways.

First write $\lambda=(\lambda_1,\ldots, \lambda_k),$ where $k=l(\lambda)$. For each row $i$ of the diagram of $\lambda$ we get $\lambda_i$ positive integers given by the entries in the $i$th row by the filling $J$. This gives us a multi-index $J_i$. \textit{The product over rows} of the filling $J$ of shape $\lambda$ is defined as 
\begin{equation}
    \Delta_{J}^r:=\prod_{1\leq i\leq k}\Delta_{J_i}.
\end{equation}
By transposing $\lambda$ we get another partition $\lambda^T$. The filling $J$ gives a filling $J^T\in \mathcal{FILL}(\lambda^T,[n])$. 
\textit{The product over columns} of the filling $J$ of shape $\lambda$ is defined as 
\begin{equation}
    \Delta^c_{J}:=\Delta^r_{J^T}
\end{equation}
\end{construction}
\begin{remark}\label{remarkaboutspan}
The classical result is that for a finte free $R$-module $E$, we can identify $\SS^{\lambda}(E)$ as the $R$-module generated by the set $\{\Delta_{J}^c\}_{J\in \mathcal{FILL}(\lambda,[n])}$, see \cite[Chapter 8]{fulton}.
\end{remark}
\begin{definition}
A semistandard Young's tableaux of shape $\lambda$ is a filling of the Young's diagram of shape $\lambda$ with positive integers so that the entries strictly increase down each column and weakly increase from left to right across each row.
\end{definition}

\begin{example}\label{exssyt}
The filling in Example \ref{wrongssyt} is not semistandard. A filling of shape $(5,3,1)$ which is semistandard is given by
$$\ytableausetup{nosmalltableaux}
\begin{ytableau}
1 & 7 & 7& 9&9 \\
2 & 8  &8\\
5\\
\end{ytableau}.$$
\end{example}
\begin{notation}\label{notationSSYT}
For a fixed $\lambda$, we denote the set of all semistandard Young's tableaux of shape $\lambda$ with entries in $[n]$ as $\mathcal{SSYT}(\lambda, [n])$.
\end{notation}

\begin{remark}
In the language of Notation \ref{notationSSYT} and Construction \ref{constfilling}, the Pl\"ucker products $(\Delta^c_J)_{J\in \mathcal{SSYT}(\lambda,[n])}$ is a basis for the finite free $R$ module $\mathbb{S}^\lambda(E)$.
\end{remark}

\subsection{Reverse Row Strict Tableaux}
For the purposes of this paper it is useful to define another filling of Young's diagrams, called the \textit{reverse row-strict tableaux}. 
\begin{definition}\label{definitionrsst}
A reverse row-strict tableau (RRST) of shape $\lambda$ is a filling of the Young's diagram of shape $\lambda$ with positive integers so that the entries weakly decrease down each column and strictly decrease from left to right along each row.
\end{definition}
\begin{example}\label{explerrst}
The filling in Example \ref{wrongssyt} is not reverse row-strict and neither is the one in Example \ref{exssyt}. A filling of shape $(5,3,1)$ which is reverse row-strict is 
$$
\ytableausetup{nosmalltableaux}
\begin{ytableau}
6 & 5 & 4& 3&1 \\
6 & 4  &2\\
1\\
\end{ytableau}.$$
\end{example}
\begin{notation}\label{reverserow}
For a fixed $\lambda$, we denote the set of all reverse row-strict tableaux of shape $\lambda$ with entries in $[n]$ as $\mathcal{RRST}(\lambda, [n])$.
\end{notation}

We now record an immediate lemma for future reference.

\begin{lemma}\label{easylemma}
Given a partition $\lambda$, transposing defines a bijection
$$\mathcal{SSYT}(\lambda,[n])\to \mathcal{RRST}(\lambda^T, [n])$$ which sends the entry $i$ to the entry $n+1-i$.
\end{lemma}
\begin{proof}
Transposing sends the row entries to the column entries and vice-versa and also reverses the ordering between the entries. Thus it is clear that the image is a reverse row-strict tableau of shape $\lambda^T$. Transposing and sending the entry $i$ to $n+1-i$ defines an inverse.
\end{proof}

\begin{remark}\label{spanrrst}
It follows from Remark \ref{remarkaboutspan}, that the polynomials $\{\Delta_J^r\}_{J\in \mathcal{RRST}(\lambda^T,[n])}$ live in the embedded image of $\SS^{\lambda}(E)$ inside $R[Z]$.
\end{remark}
\begin{notation}\label{bigrrst}
We now use Notation \ref{reverserow} to define
\begin{equation}
 \mathcal{RRST}([n]):=\bigsqcup_{\substack{ \lambda \text{ is a partition }\\ l(\lambda)\leq n}} \mathcal{RRST}(\lambda, [n])  
\end{equation}

\end{notation}

\subsection{Composition diagrams and composition tableaux}

\begin{definition}
A composition of a non-negative integer $m$ is an ordered tuple $\alpha=(\alpha_1,\ldots, \alpha_l)$ strictly positive integers  such that $\sum_i\alpha_i=m.$ The length $l(\alpha)$ of a composition $\alpha$ is the size of the tuple $\alpha=(\alpha_1,\ldots, \alpha_l)$. 
\end{definition}

\begin{example}
Two distinct composition of $4$ of length $2$ are $(1,3)$ and $(3,1)$. The number $0$ has exactly $1$ composition, the empty composition.
\end{example}

\begin{definition}
The composition diagram associated to the composition $\alpha=(\alpha_1,\ldots, \alpha_l)$ is left justified array containing $\alpha_i$ cells  in the $i$th row.
\end{definition}

\begin{example}
The composition diagram associated to $(5,3,1)$ and $(1,5,3)$ respectively are $$\ytableausetup{centertableaux}\ytableausetup{nosmalltableaux}
\ydiagram{5,3,1}  \qquad \ytableausetup{centertableaux}\ytableausetup{nosmalltableaux}
\ydiagram{1,5,3}.$$
\end{example}

\begin{definition}
Given a composition $\alpha$, a filling of shape $\alpha$ is a filling of the cells of the composition diagram with entries from the positive integers.
\end{definition}

\begin{example}
A filling of the composition $(1,5,3)$ is given by
$$\ytableausetup{nosmalltableaux}
\begin{ytableau}
9 \\
3 & 4  &1 &6 &2\\
1&5&7\\
\end{ytableau}.$$
\end{example}
\begin{notation}\label{fillrsct}
For a fixed composition $\alpha$, we denote the set of all fillings of fixed shape $\alpha$ with values in $[n]$ as $\mathcal{CFILL}(\alpha, [n])$.
\end{notation}
\begin{construction}\label{prodrsct}
Corresponding to $J\in\mathcal{CFILL}(\alpha, [n]) $ one can define a product of Pl\"ucker coordinates in two ways.
Write $\alpha=(\alpha_1,\ldots, \alpha_l)$ where $l=l(\alpha)$. Then then moving down the rows of $\alpha$ we get $\alpha_i$ positive integers which correspond to a multi-index $J_i$. Then the \textit{product over rows } of the filling $J$ of shape $\alpha$ is defined as
\begin{equation}
    \Delta_J^r:=\prod_{1\leq i\leq l}\Delta_{J_i}.
\end{equation}

Given $\alpha$, write $\alpha^c$ for the tuple $(\alpha^c_1,\ldots, \alpha^c_{p})$ where $\alpha^c_i$ counts the number of cells going down the $i$th column of the composition diagram of $\alpha$ and $p$ is set to be $\operatorname{max}(\alpha_i)$. Thus going across the columns we get $\alpha^c_i$ integers in each column corresponding to the entries of $J$. This gives a multi-index $J^c_i$. Then the  \textit{product over columns} of the filling $J$ of shape $\alpha$ is defined to be 
\begin{equation}
    \Delta_J^c:=\prod_{1\leq i\leq p}\Delta_{J^c_i}
\end{equation}
\end{construction}

\begin{note}
While we use the same notations in Constructions \ref{constfilling} and \ref{prodrsct}, for the purposes of this paper, it will be clear from context whether we mean a filling of a Young's diagram or a composition diagram.
\end{note}
\begin{definition}\label{def:1}

Given a composition $\alpha=(\alpha_1,\ldots, \alpha_l)$ with $l(\alpha)=l$ and $\operatorname{max}(\alpha_i)=w$, we define a row-strict composition
tableau (RSCT), $Y$ of shape $\alpha$, to be a filling of the cells of $\alpha$ with positive integers from $[n]$ such that
\begin{enumerate}
\item  The entries of $Y$ strictly decrease in each row when read from left to right,

\item The entries in the leftmost column of $Y$ weakly increase when read from top to bottom,
\item and $Y$ satisfies the row-strict triple rule.

Here we say that $Y$ satisfies the row-strict triple rule if when we supplement $Y$ by adding enough cells with zero-valued entries to the end of each row so that the resulting supplemented tableau, $\overline{Y}$, is of rectangular shape $w \times l$, then for $1 \leq i < j \leq w$ and $2 \leq k \leq l$, we have
$$\overline{Y}(j, k) > \overline{Y}(i, k) \Rightarrow\overline{Y}(j, k) \geq \overline{Y}(i, k - 1).$$
\end{enumerate}
\end{definition}
\begin{example}\label{examplrsct}
Since the RSCT fillings are slightly confusing, we now give three examples of shape $(1,5,3)$, $(4,2,3)$ and  $(2,4,1)$ respectively,
 $$\ytableausetup{nosmalltableaux}
\begin{ytableau}
1 \\
9 & 7  &4 & 3&2\\
10&9&7\\
\end{ytableau} \qquad
\begin{ytableau}
5& 3 & 2 & 1 \\
7& 3 \\
8& 7 & 4  
\end{ytableau} \qquad
\begin{ytableau}
4& 3 \\
7& 4 &3&1 \\
9
\end{ytableau}  .$$
\end{example}

\begin{example}
We now give non-examples of RSCT fillings which each fail the row-strict triple rule in some way. We will simply tweak the RSCT fillings in Example \ref{examplrsct}. We have
 $$\ytableausetup{nosmalltableaux}
\begin{ytableau}
1 \\
10 & 7  &4 & 3&2\\
10&9&7\\
\end{ytableau} \qquad
\begin{ytableau}
5& 3 & 2 & 1 \\
7& 6 \\
8& 4& 2  
\end{ytableau} \qquad
\begin{ytableau}
4& 3 \\
7& 4 &2&1 \\
9
\end{ytableau} .$$

In the first tableau, the second row of the first column has a $10$. The $9$ in the third row is greater than the $7$ on top of it in the second row, but is less than the entry to the left of $7$ which is a $10$. In the second tableau the $4$ in the third row is greater than the $3$ in the top row but is less than the $5$ on the left of $3.$ In the last example the $2$ in the second row is greater than the $0$ on top of it (when the tableau is completed to a rectangle), but is less than the $3$ on the left of it.
\end{example}

\begin{notation}\label{notationrsct}
We denote all RSCT fillings of shape $\alpha$ as $\mathcal{RSCT}(\alpha, [n]).$
\end{notation}

\begin{remark}
In the language of Construction \ref{prodrsct} and Notation \ref{notationrsct}, the top justified minors labelled by RSCT fillings of shape $\alpha$ are given by the set $\{\Delta^r_J\}_{J\in \mathcal{RSCT}(\alpha, [n])}.$
\end{remark}

\begin{notation}\label{zeroset}
Now note that given a composition $\alpha$, there is exactly one partition, denote $P(\alpha)$ that one can obtain by rearranging the entries of $\alpha$ in a weakly decreasing order. Given a partition $\lambda$, we denote the set of all compositions which map under this transformation to $\lambda$ as the set 
\begin{equation}
    \cZ(\lambda):=\{\alpha \ | \alpha\  \textrm{is a composition and}\ P(\alpha)=\lambda \}
\end{equation}
In the sequel it will be useful to keep track of the RSCT fillings of those $\alpha$ which lie in $\cZ(\lambda).$
Therefore we set 
\begin{equation}
    \mathcal{RSCT}(\lambda,[n]):=\bigsqcup_{\alpha\in \cZ(\lambda)}\mathcal{RSCT}(\alpha,[n]).
\end{equation}

Note that this is consistent with our notation in Theorem \ref{theorem1}.
\end{notation}
\begin{remark}\label{spanrsct}
The result that we want to prove is that the set $\{\Delta_{J}^r\}_{J\in \mathcal{RSCT}(\lambda^T,[n])}$ forms a basis for the Schur module $\SS^{\lambda}(E)$. Indeed, it follows from Remark \ref{remarkaboutspan} that these polynomials live in the correct submodule of $R[Z]$.
\end{remark}

\begin{notation}\label{bignotation}
For the sequel using Notation \ref{zeroset} we define 

\begin{equation}
    \mathcal{RSCT}([n]):=\bigsqcup_{\substack{ \lambda \text{ is a partition }\\ l(\lambda)\leq n}} \mathcal{RSCT}(\lambda, [n])
\end{equation}
\end{notation}

\subsection{Monomial term orders}\label{monotermoder}

\begin{definition}
A monomial order on a polynomial ring is an order on the monic monomials of the ring such that whenever $v\leq w$, then $vt\leq wt$, for all monomials $t$ in the polynomial ring. 
\end{definition}

\begin{remark}\label{rmkini}
Let $S=R[X_1,\ldots,X_k]$ be a polynomial ring. Then a monomial order on $S$ allows us to define a natural notion of an \emph{initial term} of polynomials in $S$.
\end{remark}
\begin{definition}
Let $S$ be as in Remark \ref{rmkini} with a fixed monomial order of $S$ and let $p(X)$ be in $S\setminus \{0\}$. Write $p(X)=\sum_{u(X)\in M}c_{u(X)} u(X)$ with $M$ being some non-empty subset of the set of monomials of $S$ with $c_{u(X)}\neq 0$.
Then the initial term, denoted $i(p(X))$, of $p(X)$ is the highest order monomial $u(X)\in M$.
\end{definition}

\begin{remark}\label{multiini}
Note that taking the initial term of a polynomial is multiplicative i.e. in the notation of Remark \ref{rmkini} if $p(X),q(X)\in S$ are two non-zero polynomials, then $i(p(X)q(X))=i(p(X))i(q(X))$
\end{remark}
\begin{remark}
In our case the ring is $R[Z]$ the polynomial ring in the $n\times n$ variables $Z_{i,j}$. Since we are focusing on the Pl\"ucker algebra it is enough for us to choose a monomial order which works well with the Pl\"ucker coordinates, and more generaly matrix minors. For this we can choose an \textit{anti-diagonal} term order.
\end{remark}

\begin{definition}\label{dfn:5}
A anti-diagonal term order on $R[Z]$ is a monomial order such that in any minor of the matrix $[Z_{i,j}]$ the anti-diagonal monomial is the initial term.
\end{definition}

\begin{example}\label{exampleontermorder}
Consider the Pl\"ucker Coordinate $\Delta_{(1,5,7)}$
We write $$\Delta_{(1,5,7)}=\operatorname{det}\begin{pmatrix}
    Z_{1,1} & Z_{1,5} &Z_{1,7}  \\
    Z_{2,1} & Z_{2,5} &Z_{2,7}\\
    Z_{3,1}&Z_{3,5}& Z_{3,7}
    \end{pmatrix}$$
Under any anti-diagonal term order the leading term is  ${Z}_{3,1}Z_{2,5}Z_{1,7}$.
\end{example}

\begin{remark}\label{rmkiniterm}
Given a anti-diagonal monomial order on $\CC[Z]$, Example \ref{exampleontermorder} shows that it is easy to read off the the initial term of Pl\"ucker coordinates. For example let $\Delta_{(i_1,\ldots, i_k)}$ be a non-zero Pl\"ucker coordinates then $i(\Delta_{(i_1,\ldots, i_k)})=\prod_{m=1}^{m=k}Z_{k-m+1,i_m}$. For general products of Pl\"ucker coordinates the initial term is the product of the initial terms of the factors. In fact as the next two examples show, the products that we shall use in this note have an easy description.
\end{remark}
\begin{example}\label{initermrrst}
Let $\lambda$ be a partition and let $\mathcal{RRST}(\lambda,[n])$ be as in Notation \ref{reverserow}. Let $J\in \mathcal{RRST}(\lambda,[n])$, then $\Delta_J^r$. Pick any anti-diagonal term order on $R[Z]$
\begin{equation}\label{equationeexplicitformrrst}
    i(\Delta_J^r)=\prod_{(i,j) \text{ is a cell of } J} Z_{j,J(i,j)}=\prod_{j} \quad \big (\prod_{\substack{k \text{ is an entry}\\ \text{ in column } j \text{ of } J}}Z_{j,k} \big).
\end{equation}
\end{example}

\begin{example}
Let $T$ be the RRST filling in Example \ref{explerrst},
$$
\ytableausetup{nosmalltableaux}
\begin{ytableau}
6 & 5 & 4& 3&1 \\
6 & 4  &2\\
1\\
\end{ytableau}.$$
Then the initial term of $\Delta_T^r$ is given by 
\begin{equation}
    i(\Delta_T^r)=(Z_{1,1}Z_{1,6}^2)(Z_{2,4}Z_{2,5})(Z_{3,2}Z_{3,4})(Z_{4,5})(Z_{5,1}).
\end{equation}

\begin{example}\label{initermrsct}
Let $\alpha$ be a composition and let $\mathcal{RRST}(\alpha,[n])$ be as in Notation \ref{notationrsct}. Let $J\in \mathcal{RCST}(\alpha,[n])$, then $\Delta_J^r$. Pick any anti-diagonal term order on $R[Z]$
\begin{equation}\label{equationexplcitformrsct}
    i(\Delta_J^r)=\prod_{(i,j) \text{ is a cell of } J} Z_{j,J(i,j)}=\prod_{j} \quad \big (\prod_{\substack{k \text{ is an entry}\\ \text{ in column } j \text{ of } J}}Z_{j,k} \big).
\end{equation}
\end{example}
\end{example}

\begin{example}
Let us consider the first RSCT filling $F$ of Example \ref{examplrsct} given by $$\ytableausetup{nosmalltableaux}
\begin{ytableau}
1 \\
9 & 7  &4 & 3&2\\
10&9&7\\
\end{ytableau}.$$ Then the initial term of $\Delta_F^r$ under any anti-diagonal term order is given by $$i(\Delta_F^r)=(Z_{1,1}Z_{1,9}Z_{1,10})(Z_{2,7},Z_{2,9})(Z_{3,4}Z_{3,7})(Z_{4,3})(Z_{5,2}).$$
\end{example}

\section{Proof of the theorem}

In this section we will fix a commutative ring $R$, a finite free $R$-module $E$ and a partition $\lambda$ of a positive integer. 

\begin{remark}
While we choose to work with an arbitrary commutative ring, for concreteness we may also work over $\ZZ$, and the proof follows formally by the base change properties of the Schur functors. However, we do not use this reduction because it doesn't seem particularly necessary.
\end{remark}

A first step in the direction of the proof is to establish a bijection between the set of RRST fillings of $\lambda^T$ (Definition \ref{definitionrsst}) and the RSCT fillings (Definition \ref{def:1}) of all compositions $\alpha$ such that $\alpha\in \cZ(\lambda^T)$ (Notation \ref{zeroset}).

For this we will use \cite[Section 3.1]{mason} to define a map 
\begin{equation}\label{naturalmap}
\rho:\mathcal{RRST}(\lambda,[n])\to \mathcal{RSCT}(\lambda,[n]),    
\end{equation}
 (see Notations \ref{definitionrsst} and \ref{notationrsct}) which turns out to be a bijection.

\begin{proposition}\label{biigproposition}
There exists a maps $\rho$ as in (\ref{naturalmap}) which induces an isomorphism between the source and the target sets.
\end{proposition}
\begin{proof}
For our purposes we show that such a map exists, following \cite[Section 3.1]{mason}. The fact that the image of this morphism is correct and it is an isomorphism is Lemma 3.4 in \textit{loc. cit.}.

Given a filling $F\in \mathcal{RRST}(\lambda,[n])$ we define a composition tableau $\rho(F)=J$ as follows
\begin{enumerate}
    \item Take the left most column of $F$ and place it in weakly increasing order (from top to bottom) in the first column of $J.$
    \item After the first $k-1$ columns have been filled place the entries of the $k$th column by beginning with the largest. Place each entry $x$ into the cell $(i,k)$ such that $(i,k)$ doesn't already contain an entry from $F$ and such that $(i,k-1)$ is strictly larger than $x.$ 
\end{enumerate}

\end{proof}

\begin{example}
We have the correspondence 
$$F=
\ytableausetup{nosmalltableaux}
\begin{ytableau}
6 & 5 & 4& 3&1 \\
6 & 4  &2\\
5\\
\end{ytableau}\mapsto \rho(F)=\begin{ytableau}
5 & 4 & 2\\
6 & 5 & 4 & 3& 1 \\
6\\
\end{ytableau},  $$
where $F\in \mathcal{RRST}(\lambda,[n]) $ and it is easy to check that $\rho(F)\in \mathcal{RSCT}(\lambda,[n])$.
\end{example}

We now record a corollary of Proposition \ref{biigproposition}. Recalling Notation \ref{bigrrst} and Notation \ref{bignotation} we have

\begin{corollary}
The map $\rho$ of Proposition \ref{biigproposition} induces a bijection
$\mathcal{RRST}([n])\cong \mathcal{RSCT}([n])$.
\end{corollary}
\begin{proof}
This is immediate from the definitons of $\mathcal{RRST}([n])$ and $\mathcal{RSCT}([n])$.
\end{proof}

Recalling Notation \ref{notationSSYT} we also have
\begin{corollary}\label{bijssytrsct}
There are bijections $\mathcal{SSYT}(\lambda,[n])\cong \mathcal{RRST}(\lambda^T,[n])\cong  \mathcal{RSCT}(\lambda^T,[n])$.
\end{corollary}

\begin{proof}
The first bijection follows from Lemma \ref{easylemma} and Proposition \ref{biigproposition}.
\end{proof}

We want to now establish that taking the Pl\"ucker products over rows of RSCT fillings (see Construction \ref{prodrsct}) is a lossless procedure i.e. the set $\{\Delta^r_F\}_{F\in \mathcal{RSCT}(\alpha,[n])}$ is in bijection with the set $\mathcal{RSCT}(\alpha,[n])$. For this we will use the notion of initial terms and monomial orderings we described in Section \ref{monotermoder}.

Recall that for a polynomial $p(Z)\in R[Z]$, we denote by $i(p(Z)$ the initial term with respect to some chosen anti-diagonal term order.

First we begin by establishing the uniqueness of the leading term for any given RRST filling. For this recall Notation \ref{bigrrst}.

\begin{lemma}\label{uniquerrst}
Let $A,B\in \mathcal{RRST}([n]).$ Then if $i(\Delta_A^r)=i(\Delta_{B}^r)$ then $A=B$.
\end{lemma}
\begin{proof}
By the description of $i(\Delta_{A}^r)$ in Example \ref{initermrrst} we see that for each row $j$ of $A$, we can rebuild the complete multi-set of entries in the $j$th column of $A$ by looking at the indeterminate of form $Z_{j,k}$ showing up in $i(\Delta_{A}^r)$. The exponent of $Z_{j,k}$ counts the number of time $k$ occurs in column $j$ of $A.$ Given this multi-set there is exactly one way to write down the column $j$ of $A$, which is by arranging the multi-set in a weakly increasing order.
Thus if $A$ and $B$ are such that $i(\Delta_{A}^r)=i(\Delta_B^r)$ then $A=B$.
\end{proof}

Now we prove a lemma which will allow us to deduce the uniqueness of the leading term for RSCT fillings.

\begin{lemma}\label{initermequal}
Fix any anti-diagonal term order on $R[Z]$.
Let $J\in \mathcal{RRST}([n])$ and let $\Delta_{J}^r$ be the product as in Construction \ref{constfilling}. Let $\rho(J)\in \mathcal{RSCT}([n])$ and let $\Delta^r_{\rho(J)}$ be as in  Construction \ref{fillrsct}. Then $i(\Delta_{J}^r)=i(\Delta_{\rho (J)}^r)$ where $\rho$ is the map of Proposition \ref{biigproposition}.
\end{lemma}
\begin{proof}

We have a good description of the initial term of $\Delta_J^r$ in Example \ref{initermrrst}, in particular a factorisation of the initial terms by columns of $J$. The exponent of $Z_{j,k}$ in $i(\Delta_{J}^r)$ counts the number of times $k$ occurs in column $j.$

By Example \ref{initermrsct}, the same description goes for the initial term of $\Delta_{\rho(J)}^r$ i.e. for each column $j$ of of $\rho(J)$ the exponent of $Z_{j,k}$ counts the number of times $k$ occurs in column $j$ of $\rho(J).$

Now note that the map in Proposition \ref{biigproposition} doesn't alter the row entries in $J$, but just rearranges them. In particular the number of columns of $J$ and $\rho(J)$ aree equal and also have the same multi-set of entries in the $j$th column.
Thus for each column $j$ of $J$ we have an equality

\begin{equation}
    \prod_{\substack{k \text{ is an entry }\\ \text{ in column } j \text{ of } J }}Z_{j,k}= \prod_{\substack{k \text{ is an entry }\\ \text{ in column } j \text{ of } \rho(J) }}Z_{j,k}
\end{equation}

Since $J$ and $\rho(J)$ have the same number of columns the proof follows from the explicit formulas (\ref{equationeexplicitformrrst}) and (\ref{equationexplcitformrsct}) in Examples \ref{initermrrst} and \ref{initermrsct} respectively.
\end{proof}

As a corollary of the above lemma we obtain

\begin{corollary}\label{uniqueforrsct}
Let $P,Q\in \mathcal{RSCT}([n])$ be such that $i(\Delta_{P}^r)=i(\Delta_Q^r)$, then $P=Q$.
\end{corollary}
\begin{proof}
Let $P$ and $Q$ be such that the hypothesis holds. Write $P=\rho(A)$ and $Q=\rho(B)$ for some $A,B\in \mathcal{RRST}(\lambda,[n])$. Then the initial terms of $i(\Delta_A^r)=i(\Delta_B^r)$ by Lemma \ref{initermequal} and so by Lemma \ref{uniquerrst} we have $A=B$. Since $\rho$ is a bijection, this then shows that $P=Q$.
\end{proof}

To proceed further we need to establish a bijection between $\mathcal{SSYT}(\lambda,[n])$ (see Notation \ref{notationSSYT}) and $\mathcal{RSCT}(\lambda,[n])$.

Now note that taking products for Pl\"ucker products over columns of SSYT fillings (see Construction \ref{constfilling}) is a lossless procedure i.e. the set $\{\Delta^c_J\}_{J\in \mathcal{SSYT}(\lambda,[n])}$ is in bijection with $\mathcal{SSYT}(\lambda,[n])$. The analogous result for RSCT fillngs follows from Corollary \ref{uniqueforrsct}. In particular using Corollary \ref{bijssytrsct} we can establish that

\begin{corollary}\label{cardinality}
The cardinality of the sets $\{\Delta^c_J\}_{J\in \mathcal{SSYT}(\lambda,[n])}$ and $\{\Delta^r_F\}_{F\in \mathcal{RSCT}(\lambda^T,[n])}$ is the same and is equal to the rank of the finite free module $\SS^\lambda(E).$
\end{corollary}

\begin{proof}
We only need to show that the rank of $\SS^\lambda(E)$ is $\{\Delta^c_J\}_{J\in \mathcal{SSYT}(\lambda,[n])}$, but this is exactly the classical result, see \cite[Chapter 8]{fulton}.
\end{proof}

Now let us prove a useful algebraic lemma to so that we may conclude the proof of Theorem \ref{theorem1}.

\begin{lemma}\label{algebraiclemma}
Let $S$ be a commutative ring and $X$ a finite set. Let $M$ be a free module of rank $|X|$. Let $\{e_x\}_{x\in X}$ be a set of elements of $M$ and let $\{f_x\}_{x\in X}$ be elements of $M^\vee=\operatorname{Hom}_S(M,S)$. Choose an ordering on $X$. If the matrix $[f_x(e_y)]_{(x,y)\in X\times X}$ is invertible, then $\{e_x\}_{x\in X}$ forms a basis for $M.$ Moreover this result is independent of the ordering chosen.
\end{lemma}
\begin{proof}
Choose a basis $\{u_x\}_{x\in X}$ for the module $M.$ Now each $e_y$ can be written as $e_y=\sum_{z\in X}c_{z,y}u_z$ for some $\{c_{y,z}\}_{(y,z)\in X\times X}\subset S$. The action of $f_{x}$ on $e_y$ then is 
\begin{equation}
    f_x(e_y)=\sum_{z\in X} c_{z,y}f_x(m_z).
\end{equation}

Thus we can write the matrix 
\begin{equation}\label{matrixfactor}
    [f_x(e_y)]_{(x,y)\in X\times X}=[\sum_{z\in X}f_x(m_z)c_{z,y}]_{(x,y)\in X\times X}=[f_x(m_y)]_{(x,y)\in X\times X}[c_{x,y}]_{(x,y)\in X\times X},
\end{equation}

where the last equality in the above equation is a factorisation into the product of the matrices $[f_x(m_y)]_{(x,y)\in X\times X}$ and $[c_{x,y}]_{(x,y)\in X\times X}.$ Now the matrix on the far left of (\ref{matrixfactor}) is invertible and so too must be each factor on the far right of (\ref{matrixfactor}), for otherwise taking the determinant we would get $0$ on the RHS but non-zero on the LHS.

Thus the change of basis matrix $[c_{x,y}]_{(x,y)\in X\times X}$ is invertible, which implies our result. This is independent of the ordering on $X$ because any other order changes the determinant of the matrix $[f_x(e_y)]_{(x,y)\in X\times X}$ by a sign, in particular if it is invertible with respect to one ordering, it is with another.
\end{proof}

Now we can finally prove Theorem \ref{theorem1}.

\begin{proof}(of Theorem \ref{theorem1}) Note that from Remark \ref{spanrsct}, the set $\{\Delta_{J}^r\}_{J\in \mathcal{RSCT}(\lambda^T,[n])}$ lives inside the submodule generated by $\{\Delta_{J}^c\}_{J\in \mathcal{SSYT}(\lambda,[n])}$. 
Order the set $\mathcal{RSCT}(\lambda^T,[n])$ so that $T\leq J$ if and only if $i(\Delta^r_T)\leq  i(\Delta_J^r)$. Let $\{f_J\}_{J\in \mathcal{RSCT}(\lambda^T,[n])}$ be elements of $\operatorname{Hom}_R(\SS^\lambda(E),R)$ so that 
\begin{equation}
    f_J(\Delta^r_{F})=\text{coefficient of }i(\Delta_F^r) \text{ if } F=J \text{ and } 0 \text{ otherwise}.
 \end{equation}

Corollary \ref{uniqueforrsct} implies that the matrix $f_J(\Delta^r_F)_{(R,F)\in \mathcal{RSCT}(\lambda^T,[n])^{\times 2}}$ is upper triangular with ones on the diagonal. Thus it is invertible. Then Lemma \ref{algebraiclemma} implies that $\{\Delta^r_F\}_{F\in \mathcal{RSCT}(\lambda,[n])}$ is a a basis for $\SS^\lambda (E)$.
\end{proof}

\begin{proof} (of Corollaries \ref{repcor} and \ref{basisplucker}) Corollary \ref{repcor} is immediate from Theorem \ref{theorem1}. Corollary \ref{basisplucker} follows because it is a direct sum of all representations $\SS^{\lambda}(E)$ where Corollary \ref{repcor} holds true.
\end{proof}

\begin{proof}(of Theorem \ref{theorem2})
Using Remark \ref{spanrrst} and Lemma \ref{uniquerrst}, the proof of Theorem \ref{theorem1} applies mutatis mutandis.
\end{proof}

\end{document}